\documentclass[reqno,11pt]{amsart}
\usepackage{amscd,amssymb}
\usepackage[hidelinks]{hyperref}

\theoremstyle{plain}
\newtheorem{Thm}[subsection]{Theorem}
\newtheorem{Cor}[subsection]{Corollary}
\newtheorem{Lem}[subsection]{Lemma}
\newtheorem{Prop}[subsection]{Proposition}
\newtheorem{Conj}[subsection]{Conjecture}

\theoremstyle{definition}
\newtheorem{Def}[subsection]{Definition}

\theoremstyle{remark}

\newtheorem{Rem}[subsection]{Remark}

\newtheorem{conj}{Conjecture}

\numberwithin{equation}{section}
\renewcommand{\rm}{\normalshape}

\newif\ifShowLabels
\ShowLabelstrue
\newdimen\theight
\def\TeXref#1{%
    \leavevmode\vadjust{\setbox0=\hbox{{\tt
        \quad\quad  {\small \rm #1}}}%
    \theight=\ht0
    \advance\theight by \lineskip
    \kern -\theight \vbox to
    \theight{\rightline{\rlap{\box0}}%
    \vss}%
    }}%

\ShowLabelsfalse

\renewcommand{\sec}[2]{\section{#2}\if\relax\detokenize{#1}\relax\else\label{S:#1}\fi%
    \ifShowLabels \TeXref{{S:#1}} \fi}
\newcommand{\ssec}[2]{\subsection{#2}\if\relax\detokenize{#1}\relax\else\label{SS:#1}\fi%
    \ifShowLabels \TeXref{{SS:#1}} \fi}

\newcommand{\sssec}[2]{\subsubsection{#2}\if\relax\detokenize{#1}\relax\else\label{SSS:#1}\fi%
    \ifShowLabels \TeXref{{SSS:#1}} \fi}

\newcommand{\refs}[1]{Section ~\ref{S:#1}}

\newcommand{\reft}[1]{Theorem ~\ref{T:#1}}
\newcommand{\refl}[1]{Lemma ~\ref{L:#1}}
\newcommand{\refp}[1]{Proposition ~\ref{P:#1}}
\newcommand{\refc}[1]{Corollary ~\ref{C:#1}}

\newcommand{\refe}[1]{\eqref{E:#1}}

\newenvironment{thm}[1]%
    { \begin{Thm} \label{T:#1}  \ifShowLabels \TeXref{T:#1} \fi }%
    { \end{Thm} }

\renewcommand{\th}[1]{\begin{thm}{#1} \sl }
\renewcommand{\eth}{\end{thm} }

\newenvironment{lemma}[1]%
    { \begin{Lem} \label{L:#1}  \ifShowLabels \TeXref{L:#1} \fi }%
    { \end{Lem} }
\newcommand{\lem}[1]{\begin{lemma}{#1} \sl}
\newcommand{\elem}{\end{lemma}}

\newenvironment{propos}[1]%
    { \begin{Prop} \label{P:#1}  \ifShowLabels \TeXref{P:#1} \fi }%
    { \end{Prop} }
\newcommand{\prop}[1]{\begin{propos}{#1}\sl }
\newcommand{\eprop}{\end{propos}}

\newenvironment{corol}[1]%
    { \begin{Cor} \label{C:#1}  \ifShowLabels \TeXref{C:#1} \fi }%
    { \end{Cor} }
\newcommand{\cor}[1]{\begin{corol}{#1} \sl }
\newcommand{\ecor}{\end{corol}}

\newenvironment{defeni}[1]%
    { \begin{Def} \label{D:#1}  \ifShowLabels \TeXref{D:#1} \fi }%
    { \end{Def} }
\newcommand{\defe}[1]{\begin{defeni}{#1} \sl }
\newcommand{\edefe}{\end{defeni}}

\newenvironment{remark}[1]%
    { \begin{Rem} \label{R:#1}  \ifShowLabels \TeXref{R:#1} \fi }%
    { \end{Rem} }
\newcommand{\rem}[1]{\begin{remark}{#1}}
\newcommand{\erem}{\end{remark}}

\newenvironment{conjec}[1]%
    { \begin{Conj} \label{Co:#1}  \ifShowLabels \TeXref{Co:#1} \fi }%
    { \end{Conj} }
\renewcommand{\conj}[1]{\begin{conjec}{#1} \sl }
\newcommand{\econj}{\end{conjec}}

\newcommand{\eq}[1]%
    { \ifShowLabels \TeXref{E:#1} \fi
       \begin{equation} \label{E:#1} }
\newcommand{\eeq}{ \end{equation} }

\newcommand{\prf}{ \begin{proof} }
\newcommand{\epr}{ \end{proof} }

\newcommand\alp{\alpha}

\newcommand\del{\delta}

\newcommand\kap{\kappa}
\newcommand\lam{\lambda}		
\newcommand\sig{\sigma}		

\newcommand\ome{\omega}		\newcommand\Ome{\Omega}

\newcommand\calD{{\mathcal{D}}}

\newcommand\calH{{\mathcal{H}}}

\newcommand\calL{{\mathcal{L}}}
\newcommand\calM{{\mathcal{M}}}

\newcommand\calO{{\mathcal{O}}}

\newcommand\calS{{\mathcal{S}}}

\newcommand\calU{{\mathcal{U}}}
\newcommand\calV{{\mathcal{V}}}
\newcommand\calW{{\mathcal{W}}}

\newcommand\calY{{\mathcal{Y}}}
\newcommand\calZ{{\mathcal{Z}}}

\newcommand\PP{\mathbb{P}}

\newcommand\GG{\mathbb{G}}

\newcommand\CC{\mathbb{C}}

	\newcommand\grg{{\mathfrak{g}}}

	\newcommand\grt{{\mathfrak{t}}}

\newcommand\sdp{\times \hskip -0.3em {\raise 0.3ex
\hbox{$\scriptscriptstyle |$}}} 

\newcommand\Ad{\operatorname{Ad}}

\newcommand\diag{\operatorname{diag}}

\newcommand\End{\operatorname{End\,}}
\newcommand\Ext{\operatorname{Ext}}

\newcommand\GL{\operatorname{GL}}

\newcommand\Hom{\operatorname {Hom}}

\newcommand\PGL{{\rm PGL}}

\newcommand\rk{\operatorname{rk}}

\newcommand\Res{\operatorname{Res}}

\newcommand\Tr{\operatorname{Tr}}

\newcommand\oF{{\overline{F}}}

\newcommand\x{\times}
\newcommand\ten{\otimes}

\newcommand{\ra}{\rangle}
\newcommand{\la}{\langle}

\newcommand\Bun{\operatorname{Bun}}
\newcommand\cusp{\operatorname{cusp}}
\renewcommand\PGL{\operatorname{PGL}}
\newcommand{\DAHA}{\mathsf{H}\!\mathsf{H}}
\newcommand\Pic{\operatorname{Pic}}
\newcommand{\bpo}{\Bun_G(\PP^1)_{0,\infty}}
\newcommand\Lie{\operatorname{Lie}}

\begin{document}
\title{Square-integrable representations and Hecke eigenvalues over a local non-archimedean field}
\author{Alexander Braverman and David Kazhdan}
\begin{abstract}
Let $G$ be a split connected semi-simple group over a local non-archimedean field $F$.
In this paper we study the action of Hecke operators on the Schwartz space of half-densities on the moduli stack $\bpo$ of $G$-bundles on $\PP^1$ endowed with a trivialization at $0$ and $\infty$. This space is a representation of $G(F)\times G(F)$ which contains the regular representation of $G(F)$. We show that this inclusion is an isomorphism on the cuspidal part. In addition we  show that every irreducible square-integrable representation $\pi$ of $G(F)$ defines a Hecke eigenvalue. As a byproduct we show that every irreducible square-integrable representation of $G(F)$ with an Iwahori fixed vector defines a $\CC$-point on the Calogero-Moser space associated to $G$. We also show that for $G=\PGL(N)$ the eigenvalues of certain "simplest non-trivial" Hecke operators are essentially equal to the values of the character of $\pi$ on regular elliptic elements. These results are parallel to some of the results of \cite{BK-finite} where analogous statements are proved in the case when $F$ is a finite field.
\end{abstract}
\maketitle

\sec{int}{Introduction}

\ssec{}{Schwartz spaces of moduli spaces}
Let $F$ be a local non-archimedean field with ring of integers $\calO_F$, residue field $k$ and separable closure $\oF$.
Let $G$ be a split connected semisimple group over $F$.
Let $C$ be a smooth projective algebraic curve over $F$. We denote by $\Bun_G(C)$ the moduli stack of principal $G$-bundles on $C$.
More generally, let $A=\sum a_i x_i$ be an effective divisor on $C$ defined over $F$. In this case we can consider the stack
$\Bun_G(C,A)$ classifying $G$-bundles endowed with a trivialization along the closed subscheme $A$.

Let now $X$ be a sufficiently nice algebraic stack (e.g. an algebraic stack which can locally be represented as a quotient of a smooth affine scheme of finite type over $F$ by a linear algebraic group). Recall that given a complex number $\kappa$ following \cite{GK-further,BKP} we can associate to $X$ a vector space $\calS_{\kappa}(X)$ of "$\kappa$-measures" on $X(F)$. In the case when $X$ is a smooth algebraic variety over $F$ this is just the space of locally constant compactly supported sections of the sheaf $|\Ome_X|^{\kappa}$; in particular, if $\kap=0$ this is the space of locally constant compactly supported functions on $X(F)$, and for $\kappa=1$ this is the space of locally constant compactly supported distributions on $X(F)$. We recall the details (in particular, we recall the generalization to stacks in \refs{schwartz}). We shall mostly be concerned with $\calS_{1/2}(\Bun_G(C,A))$. This is a representation of the group $G_A(F)$ where $G_A$ is the algebraic group of maps from $A$ to $G$ when we view $A$ as a closed subscheme of $C$.

\ssec{}{Hecke operators}
Let now $D$ be a "colored effective divisor" of the form $\sum \lam_j y_j$ where $y_j$ are points of $C$ over $\oF$ and $\lam_j$ are dominant coweights of $G$; we assume that $D$ is defined over $F$. Assume in addition that $D$ is disjoint from $A$. Then as in \cite{BK-icm} we can attach to $D$ the corresponding Hecke operator 
$$
h_D\in \End(\calS_{1/2}(\Bun_G(C,A)))\ten
\bigotimes_{x\in |D|}
\left|\operatorname{Nm}_{F(x)/F}\Ome_{C,x}\right|^{-\la\lam_x,\rho^{\vee}\ra}.
$$
Here $|D|$ is the set of closed points in the support of $D$, $F(x)$ is
the residue field of $x$, $\lam_x$ is the common coweight attached to
its geometric points, and $\rho^{\vee}$ is the half-sum of the positive
roots of $G$. The line $\Ome_{C,x}$ is a one-dimensional $F(x)$-vector
space; its norm is an $F$-line, obtained by descending the tensor product
over the geometric points above $x$. The absolute-value construction
then gives a complex density line, whose real powers are canonically
defined, including when $\la\lam_x,\rho^{\vee}\ra$ is a half-integer.
We use the operator-valued convention of \cite[Theorem 1.1 and
equations (1.4)--(1.6)]{EFK}, hence the negative exponent.
 If $D=D_1+D_2$ we have $h_D=h_{D_1}\circ h_{D_2}$ (in particular all the operators $h_D$ commute).
\ssec{P1}{Our setup} In this case we shall be concerned with the case when $C=\PP^1$ and $A=0+\infty$. We shall denote
$\calS_{1/2}(\Bun_G(C,A))$ by $\calV(G)$ in this case. We trivialize $\Ome_C$ away from $A$ by $dt/t$, where $t$ is the standard coordinate on $\PP^1$, so we can view the operators $h_D$ as genuine endomorphisms of $\calV(G)$. In addition the group $G_A$ is just $G\times G$, so $\calV(G)$
is a representation of the group $G(F)\x G(F)$; the operators $h_D$ commute with the $G(F)\x G(F)$-action.  In this paper we would like to study the space $\calV(G)$ with the above structures; this paper can be seen as an upgrade of our paper \cite{BK-finite} where analogous questions were studied in the case when $F$ is replaced by a finite field. Our first result is the following (it is parallel to Theorem 1.7 of \cite{BK-finite}):
\th{filt}
As a representation of $G(F)\times G(F)$ the space $\calV(G)$ has a filtration by dominant coweights of $G$, such that the associated graded piece corresponding to a coweight $\lam$ is isomorphic to the space of locally constant compactly supported sections of a certain explicit line bundle on $G/U_\lam\underset{M_\lam}\times G/U_\lam$, where $P_\lam$ is the parabolic subgroup defined by coweight $\lam$ in a standard way, $U_\lam$ -- its unipotent radical, $M_\lam$ -- its Levi factor.
\eth
\ssec{}{The cuspidal part}
Note that the stack $\Bun_G(\PP^1,0+\infty)$ contains an open subset isomorphic to $G$ corresponding to the trivial bundle on $C=\PP^1$. This defines a natural map $\calS(G)\to \calV(G)$ which is easily seen to be an embedding. \footnote{Here we use the fact that $\Ome_G$ has a canonical trivialization (up to scalar), so there is actually a choice of a Haar measure on $G(F)$; however, what is important for us is that we have a canonical isomorphism $\End(\calS(G))\simeq \End(\calS_{1/2}(G))$.} This embedding is $G(F)\times G(F)$-equivariant, but the subspace $\calS(G)\subset \calV(G)$ is not preserved by the Hecke operators. However, \reft{filt} immediately implies the following
\prop{cuspid}
For any representation $M$ of $G(F)\times G(F)$ let $M_{\cusp}$ denote its cuspidal part.
Then the embedding $\calS(G)_{\cusp}\hookrightarrow \calV(G)_{\cusp}$ is an isomorphism. In particular, $\calS(G)_{\cusp}$ is preserved by the Hecke operators.
\eprop
\ssec{}{Hecke eigenvalues}
For any representation $(\pi,V)$ of $G(F)$ we may consider $\End(V)$ as a representation of $G(F)\times G(F)$. This representation is not smooth and we denote by $\End'(V)$ its smooth part.
Then using \refp{cuspid} we can ask the following question. Let $(\pi,V)$ be an irreducible cuspidal representation of $G(F)$. Then
$\End'(V)$ is naturally an irreducible cuspidal $G(F)\x G(F)$-representation and $\dim\Hom_{G(F)\x G(F)}(\End'(V),\calS(G)_{\cusp})=1$ (the canonical map from $\End'(V)$ to $\calS(G)$ is given by matrix coefficients, which lands in $\calS(G)$ since $V$ is cuspidal and $G(F)$ has compact center).
Thus every $h_D$ as above defines a $G(F)\x G(F)$-endomorphism of $\End'(V)$, which has to be equal to multiplication by a scalar $\eta_D(\pi)$ by irreducibility.
\ssec{}{An example} Let now $G=\PGL(N)$. Let $E/F$ be a finite separable extension of degree $N$ and let $y\in E^*$ be an element which does not lie in any smaller extension. Then $y$ defines an elliptic conjugacy class $\Ome_y$ in $G(F)$; abusing slightly the notation we will often think about $y$ as an actual element of $\PGL(N,F)$ -- this will cause no harm as long as we deal with conjugation-invariant quantities.

It follows from the above assumption that the orbit of $y$ under  the absolute Galois group $\text{Gal}(\oF/F)$ consists of $N$ distinct  elements $y_1=y, y_2,\cdots, y_N$. Consider the colored divisor $D$ of the form $\sum\limits_{i=1}^N \ome_1 y_i$ where $\ome_1$ is the first fundamental coweight of $\PGL(N)$.
The following result is parallel to Theorem 1.6 of \cite{BK-finite}:
\th{GLN}
Let $\pi$ be an irreducible cuspidal representation of $G(F)$.
For every choice of a Haar measure $dg$ on $G(F)$ we have
$$
\eta_D(\pi)=c(y,dg)\frac{\chi_\pi(y)}{d_\pi(dg)}
$$
where
\begin{enumerate}
\item
$\chi_\pi$ denotes the character of $\pi$
\item $d_\pi(dg)$ is the formal degree of $\pi$ associated to the Haar measure $dg$
\item  $c(y,dg)$ is a certain conjugation-invariant function of $y$ (this function is explicitly computed in \refs{character}).
\end{enumerate}
\eth

\ssec{}{Remark} As was mentioned above the Hecke operator $h_D$ does not depend on the choice of $dg$. If we multiply $dg$ by a constant $C$, the formal degree $d_\pi(dg)$ gets multiplied by $C^{-1}$. It follows that the function $c(y,dg)$ also gets multiplied by $C^{-1}$.
\ssec{}{Extension to square-integrable representations} As we have already mentioned, the above results are very much parallel to the results of \cite{BK-finite}. One of the purposes of this paper is to extend the problem from cuspidal representations to square-integrable ones. This can be done in the following way.

 Note that when $(\pi,V)$ was cuspidal we used a natural embedding $\End'(V)\hookrightarrow \calS(G)_{\cusp}\hookrightarrow \calV(G)$. However, these embeddings canonically split as direct summands, so we get a natural surjective map $\calV(G)\to \End'(V)$.

It turns out that when $(\pi,V)$ is no longer cuspidal, this is the correct way to proceed. First note that
$\dim\Hom_{G(F)\times G(F)}(\calS(G),\End'(V))=\dim\Hom_{G(F)\times G(F)}(\calS(G),\End(V))=1$ for any irreducible representation $(\pi,V)$ of $G(F)$. One may ask if the corresponding map from $\calS(G)\to \End(V)$ can be extended to $\calV(G)$ and whether such an extension is unique. The following Theorem can be viewed as the main result of this paper:
\th{main}
\begin{enumerate}
\item
Let $(\pi,V)$ be an irreducible square-integrable representation of $G(F)$.
Then the restriction map
\[
\Hom_{G(F)\times G(F)}(\calV(G),\End(V))
\longrightarrow
\Hom_{G(F)\times G(F)}(\calS(G),\End(V))
\]
is an isomorphism. In particular,
$\dim\Hom_{G(F)\times G(F)}(\calV(G),\End(V))=1$.
\item
For $G=\PGL(N)$ the analog of \reft{GLN} holds for any square-integrable $\pi$.
\end{enumerate}
\eth

\noindent
{\bf Remark.} It is clear that in all of the above one can replace $\End(V)$ by $\End'(V)$ (since any morphism from a smooth representation to an abstract representation necessarily lands in the smooth part of the latter).

\reft{main} allows us to talk about the Hecke eigenvalues $\eta_D(\pi)$ for any square-integrable $\pi$ -- this is just the eigenvalue of the action of $h_D$ on the one-dimensional vector space $\Hom_{G(F)\times G(F)}(\calV(G),\End(V))$.
\ssec{}{The Iwahori case and points on the Calogero-Moser space} The following special case of \reft{main} seems to be very interesting. Namely, let $I$ denote an Iwahori subgroup of $G(F)$ and let $\calH(G,I)$ denote the corresponding Iwahori-Hecke algebra. Assume that $(\pi,V)$ above is such that $V^I\neq 0$.

On the other hand let $\DAHA(G)$ denote the corresponding double affine Hecke algebra "at critical level". In fact, our notations are different from standard DAHA notations: usually $\DAHA(G)$ depends on two parameters $q$ and $t$. Our algebra corresponds to $q=1$ (this is what "critical level" means); on the other hand, $t$ should be equal to the number of elements in the residue field of $F$, which we actually denote by $q$ (thus it has nothing to do with the $q$ from the standard DAHA notation -- we hope that it will not lead to a confusion, since the standard DAHA notation will never appear in this paper). The algebra $\DAHA(G)$ contains $\calH(G,I)$ as a subalgebra. Then it follows from \cite{GK-further} (Theorem 7.9) that

\medskip
(a) $\calV(G)^{I\times I}$ is a $\DAHA(G)\otimes \DAHA(G)$-module. Moreover, the action of the Hecke operators $h_D$ (recall that $D$ is assumed to be disjoint from $0$ and $\infty$) commutes with the $\DAHA(G)\otimes \DAHA(G)$-action.

(b) There is an isomorphism of vector spaces $\DAHA(G)\simeq \calV(G)^{I\times I}$ under which the left multiplication goes to the first action of $\DAHA(G)$ on $\calV(G)^{I\times I}$ and the right multiplication goes over to the 2nd action up to certain anti-involution $\kappa$ on $\DAHA(G)$. The anti-involution $\kappa$ preserves the subalgebra $\calH(G,I)$ and its restriction to $\calH(G,I)$ is equal to the anti-involution induced by the map $g\mapsto g^{-1}$ on $G$.

\medskip
\noindent
As a corollary of (a) and (b) we get that $\Hom_{\calH(G,I)-bimod}(\DAHA(G),\End(V^I))$ has dimension 1 (the Hom is taken in the category of $\calH(G,I)$-bimodules). Clearly the center $Z(\DAHA(G))$ of $\DAHA(G)$ acts on the above Hom, and since the latter has dimension 1, we see that $\pi$ defines a homomorphism $\psi_{\pi}:Z(\DAHA(G))\to \CC$.
On the other hand, according to \cite{Ob} the algebra $Z(\DAHA(G))$ is the algebra of functions on a certain normal, irreducible Cohen-Macaulay algebraic variety $\calM(G)$ of dimension $2\rk(G)$. This variety is called the Calogero-Moser space associated to $G$ (when $G$ is classical the variety $\calM(G)$ has an explicit quiver-type description). Thus $\psi_\pi$ can be regarded as a point on $\calM(G)$, i.e. we can associate a point on $\calM(G)$ to any square-integrable representation of $G(F)$ with an $I$-fixed vector. It would be very interesting to produce a recipe for computing this point in the cases when $\calM(G)$ has an explicit description (note that square-integrable representations of $G(F)$ with an $I$-fixed vector were classified in \cite{KL}).

We conclude the introduction with the following:
\conj{center-sqint}
The point $\psi_\pi$ depends only on the central character of $V^I$ considered as a $\calH(G,I)$-module.
\econj
\ssec{}{Acknowledgments} We thank R.~Bezrukavnikov and M.~Finkelberg for very helpful discussions. The research of A.B. was partially supported by NSERC. The research of D.K.  was partially supported by the
ERC grant no. 101142781.
\sec{schwartz}{Digression on Schwartz spaces for stacks}

\ssec{}{Smooth sections of line bundles on algebraic varieties}
Let $X$ be an algebraic variety over $F$. We denote by $C^{\infty}(X)$ the space of locally constant functions on $X(F)$ and by $\calS(X)$ its subspace of functions with compact support.

For a line bundle  $\calL$ over $X$ we denote by  $\calL^0:=\calL\setminus X$ the corresponding $\GG_m$-torsor over $X$ and set
$$
|\calL|^{\kap}=\calL^0(F)\underset{F^*}\times \CC_{\kap},
$$
where $\CC_{\kap}$ denotes the 1-dimensional space $\CC$ on which $F^*$ acts by $|\cdot|^{\kap}$.
Then $|\calL|^{\kap}$ is a complex line bundle over $X(F)$. Since the bundle $ |\calL|^{\kap} $ is
 locally trivial with respect to the natural topology  we can define  its space of smooth sections which we denote by $C^{\infty}(X,|\calL|^{\kap})$. Similarly,
we denote by $\calS(X,|\calL|^{\kap})\subset C^{\infty}(X,|\calL|^{\kap}) $ the subspace of  sections with compact support.

In the case when $X$ is smooth we shall often take $\calL=\ome_X$, where $\ome_X$ is the line bundle of differential forms of top degree and write $\calS_{\kap}(X)$ instead of $\calS(X,|\ome_X|^{\kap})$.
The case $\kap =1/2$ is of special interest since  the space $\calS_{1/2}(X)$ is endowed with a natural Hermitian product. We denote by $L^2(X)$ its Hilbert space completion.
\ssec{}{Remark}
\begin{enumerate}
\item If $U\subset X$ is an open subset and $Z=X\backslash U$ then we have a short exact sequence
$$
0\to \calS(U,|\calL|_U^{\kap})\to \calS(X,|\calL|^{\kap})\to \calS(Z,|\calL|_Z^{\kap})\to 0.
$$
\item
More generally, instead of choosing $\calL\in \Pic(X)$ and $\kap \in \CC$ we can start with any element of $\Pic(X)\otimes \CC$ -- all the above definitions make sense in this context.
\end{enumerate}

\ssec{}{Smooth sections on  stacks} In this
subsection we
extend the above definitions to a class of algebraic stacks.

\defe{adm}
An algebraic stack $\calY$ is {\it admissible} if it is a union of
increasing open substacks of finite type, each of which has a presentation
as a quotient stack $X/H$, where $X$ is a smooth variety and $H$ is a
smooth affine algebraic group.
\edefe

A  presentation of $\calY$ as a quotient $\calY= X/\GL(n)$
is called an {\it admissible presentation}. We denote its projection by
$p:X\to\calY$.\footnote{We use the quotient-stack convention of
\cite[Section 2.2]{BKP}, which is slightly stronger than the notion in
\cite[Section 6]{GK-further}.}

\ssec{}{Remark}\begin{enumerate}
\item Any smooth  admissible stack of finite type can be presented as a quotient  $X/\GL(n) $ for a smooth variety $X$
(see \cite[Section 2.2]{BKP}). By Hilbert's theorem 90, the groupoid
$\calY(F)$ is equivalent to the quotient groupoid $X(F)/\GL(n,F)$.
\item Any admissible stack is automatically locally of finite type.
\item A line bundle on a quotient  $X/G$  is a $G$-equivariant line bundle on $X$.

\end{enumerate}

\defe{}
\begin{enumerate}
\item
Assume that $\calY$ is an admissible stack of finite type over $F$. Choose an admissible presentation $\calY=X/\GL(n)$ for some variety $X$ and set
$$
\calS(\calY,|\calL|^{\kap})
=\calS(X,|p^*\calL|^{\kap}\otimes|\ome_{X/\calY}|)_{\GL(n,F)},
$$
where the subscript denotes coinvariants for the natural $\GL(n,F)$-action.
The relative-density factor is necessary for independence of the presentation.
\item
If  $\calY$ is only locally of finite type, then we can write $\calY$ as a direct limit of open sub-stacks $\calY_i$ of finite type over $F$
and define $\calS(\calY,|\calL|^{\kap}):=\underset{\to}\lim\  \calS(\calY_i,|\calL|^{\kap})$.
\end{enumerate}
\edefe
The above definition makes sense because we claim that if  $\calY$ is an admissible stack of finite type then the space $\calS(\calY,|\calL|^{\kap})$ does not depend on the choice of an admissible  presentation $\calY=X/\GL(n)$.
For $\calL=\ome_{\calY}$ and $\kap=1$ this claim is proved in
\cite[Section 6]{GK-further}; see also \cite[Definition 2.5 and Claim 2.6]{BKP}
for the version with a line bundle.
\ssec{}{Functoriality}
If $\calY$ is an admissible stack and $\calU$ is an open substack we have a natural map $\calS(\calU,|\calL|^{\kap})\to \calS(\calY,|\calL|^{\kap})$ which
 need not be injective in general. Similarly, if $\calZ$ is the complement of $\calU$ we have an exact sequence
\eq{rightex}
\calS(\calU,|\calL|^{\kap})\to \calS(\calY,|\calL|^{\kap})\to \calS(\calZ,|\calL|^{\kap})\to 0
\end{equation}
\refe{rightex} becomes a short exact sequence in the case when all points of $\calZ$ have unipotent stabilizers.

More generally, let $f:\calW\to \calY$ be a smooth representable map of admissible stacks and $\ome_{\calW/\calY}$ be the relative canonical bundle.
Then we have a natural ("integration over the fibers") map
\eq{int}
\iota_f:\calS(\calW,|f^*\calL|^{\kap}\otimes|\ome_{\calW/\calY}|)
\to \calS(\calY,|\calL|^{\kap}).
\end{equation}

\sec{bung}{The case of \texorpdfstring{$\bpo$}{the framed-bundle stack}}
\ssec{}{The stack \texorpdfstring{$\bpo$}{of framed bundles} and its stratification} Let $\bpo$ denote the moduli stack of $G$-bundles on $\PP^1$ with trivializations at $0$ and $\infty$. We shall be interested in the space $\calS_{\kap}(\bpo)$. Our main interest will be in the case $\kap=\frac{1}{2}$; in this case we denote the above space by $\calV(G)$. Note that $\bpo$ has a natural action of the group $G\times G$ (by changing the two trivializations), so for every $\kap$ the space $\calS_{\kap}(\bpo)$ is a smooth representation of $G(F)\times G(F)$. Also the stack $\bpo$ contains the group $G$ as the open subset corresponding to the trivial bundle. The canonical bundle of $\bpo$ is naturally trivialized when restricted to this open subset,  thus we get a natural map $\calS(G)\to \calS_{\kap}(\bpo)$. This map is injective, since the stabilizers of all points of $\bpo$ are unipotent.

To understand better the structure of $\calS_{\kap}(\bpo)$ we need to introduce the following notation. Recall that isomorphism classes of $G$-bundles on $\PP^1$ admitting a trivialization at $0$ are in one-to-one correspondence with dominant coweights of $G$; for each such coweight $\lam$ we denote by $\bpo^{\lam}$ the locally closed substack of $\bpo$ corresponding to $\lam$.

Let us choose a Borel subgroup $B$ of $G$. Then to each  $\lam$ as above we can associate a parabolic subgroup $P_{\lam}$ containing $B$. This subgroup is determined by the set of simple roots $\alp$ such that $\la \lam,\alp\ra=0$ (thus $P_\lam=G$ if $\lam=0$ and $P_{\lam}=B$ if $\lam$ is regular). We denote by $M_\lam$ the Levi group of $P_\lam$ and by $U_\lam$ its unipotent radical. Set
$$
X_\lam=G/U_{\lam}\underset{M_{\lam}}\times G/U_{\lam}.
$$
\th{global}
There exists a $G\times G$-equivariant unipotent group-scheme $W_{\lam}$ over $ X_\lam$ such that
    $\bpo^{\lam}$ is equivalent to some gerbe over $W_{\lam}$. In particular, there is map $\bpo^\lam\to X_\lam$ whose fibers are
    classifying stacks of the fibers of $W_\lam$.
\eth

This theorem is proven in \cite{BK-finite} (cf. Theorem 2.1 in {\em loc. cit.}).  For every $\lam$ as above we denote by $\bpo^{\leq \lam}$ the union of all $\bpo^\mu$ with $\mu\leq \lam$; this is an open substack of $\bpo$. One can similarly define the open substack $\bpo^{<\lam}$. Note that $\bpo^{\leq 0}=G$.

Since $\bpo^{\lam}$ is a $G\times G$-equivariant gerbe over $X_\lam$ it follows that any $G\times G$-equivariant line bundle on $\bpo^\lam$ canonically descends to $X_\lam$. Let $\calD_\lam$ be the determinant
of the vector bundle of Lie algebras of $W_\lam$.  Define also $\calL_\lam$ to be
the above descent of the pullback of $\ome_{\bpo}$ along the embedding $\bpo^\lam\hookrightarrow\bpo$.

Then \reft{global} immediately implies the following
\cor{filt}
\begin{enumerate}
\item For every $\kap\in\CC$ the spaces $\calS_\kap(\bpo^{\leq \lam})$ define a $G(F)\times G(F)$-invariant filtration of $\calS_\kap(\bpo)$  by dominant coweights. For every $\lam$ the corresponding associated graded piece is isomorphic as a $G(F)\times G(F)$-representation to
$\calS(X_\lam,|\calL_\lam|^\kap\otimes |\calD_\lam|^{-1})$.
\item
The embedding $\calS(G)\hookrightarrow \calS_\kap(\bpo)$ is an isomorphism on the cuspidal part.
\end{enumerate}
\ecor
\prf
The exact sequence \refe{rightex} and the unipotence of the stabilizers
give the first assertion. The second follows because, for $\lam\neq0$,
$P_\lam$ is proper and the representation attached to the stratum is
parabolically induced, so its cuspidal part is zero.
\epr
\ssec{}{The computation of \texorpdfstring{$\calL_\lam$ and $\calD_\lam$}{the line bundles}} In order to get further information about the representation $\calS_\kap(\bpo)$ beyond the cuspidal part, we need to compute the line bundles $\calL_\lam$ and $\calD_\lam$. Since $G\times G$ acts transitively on $X_\lam$ with stabilizer $M_\lam\ltimes (U_\lam\times U_\lam)$ every $G\times G$-equivariant line bundle $\calL$ corresponds to a character $\psi(\calL):M_\lam\to \GG_m$. We need to compute the characters $\psi(\calL_\lam)$ and $\psi(\calD_\lam)$.

For each $\lam$ as above let $\del_\lam$ be the determinant of the action of $M_\lam$ on the Lie algebra of $U_\lam$. Note that $\del_\lam$ only depends on the "regularity" of $\lam$; also the corresponding $G$-equivariant line bundle on $G/U_\lam$ is the canonical bundle. Similarly let $\phi_\lam$ be the character of $M_\lam$ whose restriction to a maximal torus $T$ of $B$ is equal to
$$
\prod\limits_{\alp\in \Phi^+} \alp^{\la\lam,\alp\ra},
$$
where $\Phi^+$ denotes the corresponding set of positive roots.
\th{bundles}
We have
$$
\psi(\calL_\lam) =\phi_\lam^2\ \ \ \ \  \psi(\calD_\lam)=\phi_\lam\del_\lam^{-1}
$$
\eth
\prf
Write $M=M_\lambda$ and $U=U_\lambda$. The stabilizer of the
standard point of $X_\lambda$ is
\[
H_\lambda=(U\times U)\rtimes\operatorname{diag}M.
\]
Since $U\times U$ has no nontrivial algebraic characters,
characters of $H_\lambda$ identify with characters of $M$.
Thus it suffices to compute the $M$-characters on the fibers
of $\calL_\lambda$ and $\calD_\lambda$ at this point. Here
$\calL_\lambda$ denotes the descent of the restriction of the
canonical bundle of $\operatorname{Bun}_G(\mathbb P^1)_{0,\infty}$
to the stratum indexed by $\lambda$.

Consider the grading defined by $\lambda$:
\[
\mathfrak g=\bigoplus_{n\in\mathbb Z}\mathfrak g_n,
\qquad
\mathfrak g_n=\{v:\operatorname{Ad}(\lambda(a))v=a^nv\}.
\]
Each $\mathfrak g_n$ is an $M$-module, and
\[
\mathfrak g_0=\mathfrak m,
\qquad
\mathfrak u=\bigoplus_{n>0}\mathfrak g_n.
\]
Let $d_n=\det(\mathfrak g_n)$, regarded as a character of $M$.
Then
\[
\delta_\lambda=\prod_{n>0}d_n,
\qquad
\varphi_\lambda=\prod_{n>0}d_n^{\,n}.
\]
The second identity follows by restriction to $T$ and also
shows that the displayed formula for $\varphi_\lambda|_T$
extends to a character of $M$.

Let $E_\lambda$ be the bundle representing the stratum, with
the standard framings at $0$ and $\infty$. With the convention
that the positive-degree root spaces correspond to $U$, its
adjoint bundle has the $M$-equivariant decomposition
\[
\operatorname{ad}(E_\lambda)
\simeq\bigoplus_{n\in\mathbb Z}
\mathfrak g_n\otimes\mathcal O_{\mathbb P^1}(n).
\]
Put
\[
A_\lambda=\operatorname{ad}(E_\lambda)(-0-\infty).
\]
The Lie algebra of the automorphism group preserving the
two framings is $H^0(\mathbb P^1,A_\lambda)$. Consequently,
the fiber of $\calD_\lambda$ has character
\[
\psi(\calD_\lambda)
=\prod_{n>0}d_n^{\,n-1}
=\varphi_\lambda\delta_\lambda^{-1},
\]
because $h^0(\mathcal O(n-2))=n-1$ for $n\geq1$, and
$h^0(\mathcal O(n-2))=0$ for $n\leq0$.

The tangent complex of the framed-bundle stack at
$E_\lambda$ is $R\Gamma(\mathbb P^1,A_\lambda)[1]$.
Its canonical-line fiber is therefore
\[
\det H^0(\mathbb P^1,A_\lambda)
\otimes\det H^1(\mathbb P^1,A_\lambda)^{-1}.
\]
Since $\chi(\mathcal O(n-2))=n-1$, this yields
\[
\psi(\calL_\lambda)=\prod_{n\in\mathbb Z}d_n^{\,n-1}.
\]
The adjoint representation of $M$ has trivial determinant,
so $d_0=1$. Moreover, $d_{-n}=d_n^{-1}$: their restrictions
to $T$ are products of opposite roots. Pairing the factors
indexed by $n$ and $-n$, we obtain
\[
\psi(\calL_\lambda)
=\prod_{n>0}d_n^{\,n-1}d_{-n}^{-n-1}
=\prod_{n>0}d_n^{\,2n}
=\varphi_\lambda^2.
\]
This proves both formulas.
\epr
\ssec{}{Critical level} For general $\kap$, the representations
\[
 \calS(X_\lam,|\calL_\lam|^\kap\otimes|\calD_\lam|^{-1})
\]
depend on $\lam$, rather than only on $P_\lam$ (for example,
$P_\lam=B$ for every regular $\lam$). The case $\kap=\frac12$ is special.
Set
$$
\begin{aligned}
\calV(G)_{\leq\lam}&=\calS_{1/2}(\bpo^{\leq\lam}),\\
\calV(G)_\lam&=\calS(X_\lam,|\calL_\lam^{1/2}\otimes\calD_\lam^{-1}|).
\end{aligned}
$$
(note that by \reft{bundles} the line bundle $\calL_\lam$ has a natural square root).
\cor{1/2}
$\calV(G)_\lam=\calS_{1/2}(X_\lam)$.
\ecor
\prf
According to \reft{bundles} the line bundle $\calL_\lam^{1/2}\ten \calD_\lam^{-1}$ corresponds to the character $\del_\lam$ which is a square root of $\ome_{X_\lam}$.
\epr

\sec{proof}{Proof of \reft{main}(1)}

\ssec{}{\reft{main}(1) and \texorpdfstring{Ext}{Ext}-vanishing} \refc{1/2} implies that in order to prove the 1st assertion of \reft{main} it is enough to show the following:

(a) For every $\lam\neq 0$ and every irreducible square-integrable representation $(\pi,V)$ of $G(F)$ we have
\eq{hom0}
\Hom_{G(F)\times G(F)}(\calS_{1/2}(X_\lam),\End'(V))=0.
\end{equation}

(b) For every $\lam\neq 0$ and every irreducible square-integrable representation $(\pi,V)$ of $G(F)$ we have
\eq{ext1}
\Ext^1_{G(F)\times G(F)}(\calS_{1/2}(X_\lam),\End'(V))=0.
\end{equation}

In fact we are going to prove a stronger result. Note that $\End'(V)$ is a smooth irreducible representation of $G(F)\times G(F)$, so \reft{exti} below is indeed a generalization of (a) and (b).

\th{exti}
Let $\lambda\neq 0$, and let $\Pi$ be an irreducible
square-integrable representation of $G(F)\times G(F)$. Then
\[
\operatorname{Ext}^{i}_{G(F)\times G(F)}
\bigl(\calS_{1/2}(X_\lambda),\Pi\bigr)=0
\qquad\text{for all }i\geq 0.
\]
\eth

\begin{proof}
Write $P=P_\lambda=MU$, and let $\overline P=M\overline U$
be the opposite parabolic. All induction and Jacquet
functors below are normalized. We take Ext groups in the
category of smooth complex representations of the indicated groups of $F$-points.

First, we have
\[
\calS_{1/2}(X_\lambda)
\simeq i_{P(F)\times P(F)}^{G(F)\times G(F)}\calS(M),
\]
where $\calS(M)$ carries the $M(F)\times M(F)$-action
\[
((m_1,m_2)\cdot f)(m)=f(m_1^{-1}mm_2).
\]
Indeed, the projection
\[
X_\lambda(F)\longrightarrow G(F)/P(F)\times G(F)/P(F)
\]
has fiber $(M(F)\times M(F))/\operatorname{diag}M(F)\simeq M(F)$.
The half-density character on the stabilizer
\[
(U(F)\times U(F))\rtimes\operatorname{diag}M(F)
\]
is $\Delta_P$ on the diagonal Levi, where
\[
\Delta_P(m)=|\del_P(m)|_F,\qquad
\del_P(m)=\det\bigl(\operatorname{Ad}(m)|_{\mathfrak u}\bigr),\quad m\in M(F).
\]
This is precisely the restriction of
$\Delta_P^{1/2}\boxtimes\Delta_P^{1/2}$, the character
normalizing induction from $P(F)\times P(F)$.

By second adjointness, $i_{P(F)\times P(F)}^{G(F)\times G(F)}$ is left
adjoint to $r_{\overline P(F)\times\overline P(F)}^{G(F)\times G(F)}$.
Both functors are exact, so induction preserves projectives
and the adjunction extends to all Ext groups. Thus, setting
\[
R=r_{\overline P(F)\times\overline P(F)}^{G(F)\times G(F)}\Pi,
\]
we obtain
\[
\operatorname{Ext}^{i}_{G(F)\times G(F)}
\bigl(\calS_{1/2}(X_\lambda),\Pi\bigr)
\simeq
\operatorname{Ext}^{i}_{M(F)\times M(F)}\bigl(\calS(M),R\bigr).
\]

Choose a uniformizer $\varpi$ of $F$, and set
\[
a=\lambda(\varpi^{-1})\in Z(M)(F),
\qquad
z=(a,a)\in Z(M\times M)(F).
\]
For every root $\beta$ occurring in
$\operatorname{Lie}(\overline U)$, we have $|\beta(a)|<1$.
Hence $z$ lies in the strictly contracting chamber for
$\overline P(F)\times\overline P(F)$. Since $G$ is semisimple and
$\lambda\neq 0$, its image modulo the center of $G(F)\times G(F)$
does not belong to a compact subgroup.

Casselman's square-integrability criterion implies that
every generalized central exponent $\xi$ of $R$ satisfies
\[
|\xi(z)|<1.
\]
Since $R$ has finite length, its generalized
central-character decomposition shows that $z-1$ acts
invertibly on $R$. On the other hand, $z-1$ acts as zero
on $\calS(M)$, because
\[
(z\cdot f)(m)=f(a^{-1}ma)=f(m).
\]

The actions of a central element on an Ext group induced
from its two arguments coincide. Consequently, the
endomorphism induced by $z-1$ on
\[
\operatorname{Ext}^{i}_{M(F)\times M(F)}(\calS(M),R)
\]
is simultaneously zero, by its action on the first
argument, and invertible, by its action on the second.
This Ext group therefore vanishes for every $i\geq 0$,
proving the theorem.
\end{proof}

\sec{character}{Proof of \reft{main}(2)}
Everywhere in this Section we have $G=\PGL(N)$.

\ssec{}{Proof of \reft{GLN}} Although \reft{GLN} is a special case of \reft{main}(2) which we prove at the end of this Section, we think that it is instructive to first write the proof of \reft{GLN} since it is considerably easier and essentially explains what is going on. The proof is a variation of the proof of Theorem 1.6(ii) from \cite{BK-finite} (where a similar result is proven when $F$ is a finite field).  Recall that our divisor $D$ comes from a field extension $E/F$ of degree $N$ and an element $y\in E^*$ which does not lie in any smaller field. To such a $y$ there corresponds an elliptic element in $G(F)=\PGL(N,F)$ defined uniquely up to conjugacy; we denote by $\Ome_y$ the corresponding elliptic orbit.
Then repeating word-by-word the proof of Lemma 5.1 from \cite{BK-finite} we get the following:
\lem{conv}There exists a $G(F)$-invariant measure $\mu_y$ on $\Ome_y$ such that the action of the operator $h_D$ on $\calS(G)_{\cusp}$ is equal to the convolution with $\mu_y$.
\elem
We are going to compute this measure later, but first let us explain why \refl{conv} implies \reft{GLN} (note that in the formulation of \reft{GLN} we don't say anything about the function $c(y,dg)$). Let us choose a Haar measure $dg$ on $G(F)$.
Choosing an element in $\Ome_y$ we get a map $p_y:G(F)\to \Ome_y$ with compact fibers. Thus we can consider the measure $(p_y)_*(dg)$. Since $p_y$ is defined uniquely up to shifts, this measure is canonically defined after the choice of $dg$.
\lem{formal}The measure $(p_y)_*(dg)$ acts in an irreducible cuspidal representation $(\pi,V)$ of $G(F)$ by the scalar
$\frac{\chi_\pi(y)}{d_\pi(dg)}$.
\elem
The proof of the lemma is straightforward (it is essentially the definition of $d_\pi(dg)$). Since any two invariant measures on $\Ome_y$ are proportional, we can write $\mu_y=c(y,dg) (p_y)_*(dg)$, which in view of \refl{conv} and \refl{formal} implies \reft{GLN}.
\ssec{}{The function \texorpdfstring{$c(y,dg)$}{c(y,dg)}}We want to finish the discussion of the cuspidal part with the computation of $c(y,dg)$. Let us introduce some notation.
First, we denote by $T$ the maximal torus in $G$ whose group of $F$-points is $E^*/F^*$.
Let also  $e=e(E/F)$ be the ramification index of the extension $E/F$ and $f=f(E/F)$ be the degree of the residue field extension (so we have $ef=N$). We denote by $K=G(\calO_F)$ the standard maximal open compact subgroup of $G(F)$. Set
$$
\calD_G(y)=\det(1-\Ad(y)| \grg/\grt),
$$
where $\grg=\Lie(G), \grt=\Lie(T)$. Lastly set
$$
\delta_E=|\operatorname{disc}(\calO_E/\calO_F)|_F,\qquad
 A_E=\frac{\prod_{j=1}^{N}(1-q^{-j})}
 {e\,\delta_E^{1/2}(1-q^{-f})}.
$$
We can now formulate
\th{cy}
We have
$$
c(y,dg)=\frac{A_E\,|\mathcal D_G(y)|_F^{1/2}}{dg(G(\calO_F))}
$$
\eth
The proof of \reft{cy} is given in the next Section. Let us just make some comments about the above formula. First note that as expected if we multiply $dg$ by a constant the function $c(y,dg)$ gets multiplied by the inverse constant.
Second assume that $E/F$ is unramified. In this case we have
$$
 A_E=\prod_{j=1}^{N-1}(1-q^{-j}).
$$
If in addition $y\in\calO_E^*$ has residue generating the residue
extension, then $|\mathcal D_G(y)|_F=1$, and
$$
 c(y,dg)=
 \frac{\prod_{j=1}^{N-1}(1-q^{-j})}{dg(G(\calO_F))}.
$$
Assume also that $\pi$ has depth $0$. That means that $\pi$ is the compact induction from $G(\calO_F)$ to $G(F)$ of an irreducible cuspidal representation $\sig$ of $G(k)$ (regarded as a representation of $G(\calO_F)$). It is known that in this case we have
$$
d_\pi(dg)=\frac{\dim(\sig)}{dg(G(\calO_F))}.
$$
Set $r=\frac{N(N-1)}{2}$. On the other hand, it is known that for any irreducible cuspidal representation $\sig$ of $G(k)=\PGL(N,k)$ we have
$$
\dim(\sig)=q^{r}\prod_{j=1}^{N-1}(1-q^{-j}).
$$
This shows that if we choose $dg$ so that $dg(G(\calO_F))=1$ then for $\pi$ as above we have $\eta_D(\pi)=q^{-r}\chi_\pi(y)$. This is similar to the statement of Theorem 1.6 in \cite{BK-finite}.
\footnote{The factor $q^{-r}$ is obtained with the normalization
$dx(\calO_F)=1$ for the measure associated to a differential form. If we change it to $dx(\calO_F)=q$ (which is more natural from a certain point of view), this additional factor disappears.}
Note, however, that essentially in all other cases we do need to correct the above formula.

\ssec{}{Extension to square-integrable representations}
For a cuspidal representation $(\pi,V)$ of $G(F)$, the preceding
argument used the embedding $\End'(V)\hookrightarrow\calS(G)$ given
by matrix coefficients. For a general square-integrable representation,
matrix coefficients need not have compact support. We instead compare
the Hecke action with orbital integration after passing to the
square-integrable quotient supplied by \reft{main}(1). Choose a Haar
measure $dg$ on $G(F)$ to identify $\calS(G)$ with the Hecke algebra.
For an
irreducible square-integrable representation $(\pi,V)$, let
\[
 \Phi_\pi:\calV(G)\longrightarrow\End'(V)
\]
be the unique $G(F)\times G(F)$-equivariant extension of
$f\mapsto\pi(f)$ supplied by \reft{main}(1).

We also use the operator $\pi(\mu_y):V\to V$ defined by orbital
integration. This can be defined in finite-dimensional spaces: if
$K\subset G(F)$ is compact open, put $E_K=\pi(e_K)$, where $e_K$
is the normalized averaging idempotent. For $v\in V^K$, set
\[
 \pi(\mu_y)v=\int_{\Ome_y}E_K\pi(x)v\,d\mu_y(x)\in V^K.
\]
By admissibility and absolute convergence of orbital integrals of
square-integrable matrix coefficients, this integral converges
coordinate by coordinate in the finite-dimensional space $V^K$.
Conjugation invariance of $\mu_y$ makes these definitions compatible
as $K$ varies and gives a $G(F)$-equivariant operator on $V$.
Its finite-dimensional blocks satisfy
\[
 E_{K_1}\pi(\mu_y)E_{K_2}
 =\int_{\Ome_y}E_{K_1}\pi(x)E_{K_2}\,d\mu_y(x).
\]

\prop{orbital-comparison}
Let $D$, $y$ and $\mu_y$ be as above, and let $(\pi,V)$ be an
irreducible square-integrable representation of $G(F)$.
For every $f\in\calS(G)$, one has the identity
of operators on $V$
\eq{orbital-comparison}
 \Phi_\pi(h_D(f))=\pi(\mu_y)\pi(f).
\end{equation}
\eprop
\prf
Let $j:G\hookrightarrow\bpo$ be the open
stratum of trivial bundles. Let $\calZ_D$ be the Hecke correspondence,
with input and output maps $p_1,p_2:\calZ_D\to\bpo$, and set
\[
 Z=p_1^{-1}(G),\qquad p=p_1|_Z:Z\to G,\qquad
 q=p_2|_Z:Z\to\bpo.
\]
Since the modification is minuscule, $p$ is proper and $q$ is smooth;
the canonical-line identification defining $h_D$ gives, after the
fixed trivialization of the line attached to $D$,
\eq{hecke-relative-line}
 |p^*\ome_G|^{1/2}
 \simeq |q^*\ome_{\bpo}|^{1/2}\otimes|\ome_{Z/\bpo}|.
\end{equation}
Thus $h_D(f)=q_!p^*f$, where \refe{hecke-relative-line} is understood
in this formula; see \cite[Section 5.1]{BKP}.

Let $\calU=q(Z)$ and $Z^\circ=q^{-1}(G)$. The stack $\calU$ is an
invariant open substack of finite type containing $G$, and
$q:Z\to\calU$ is a smooth cover. This cover is also surjective on
$F$-points up to isomorphism. Indeed, in the present $\PGL(N)$
situation, its fiber over a framed bundle is a nonempty open subset
of the space of inverse modifications
$\operatorname{Res}_{E/F}\PP^{N-1}_E$. The latter is an $F$-rational
variety with Zariski-dense $F$-points. The open subset $Z^\circ$ is
dense in $Z$.

For compact open subgroups $K_1,K_2\subset G(F)$, put
\[
 A_{K_1,K_2}(g)=E_{K_1}\pi(g)E_{K_2}
 \in\Hom_{\CC}(V^{K_2},V^{K_1}).
\]
Both fixed-vector spaces are finite-dimensional. On $G$, integration
against $A_{K_1,K_2}(g)|dg|^{1/2}$ gives the corresponding block of
the map $f\mapsto\pi(f)$. We pull this operator-valued half-density
back to $Z^\circ$ and show that it is locally integrable on all of
$Z$, coordinate by coordinate.

The analytic input is that matrix coefficients of a square-integrable
representation belong to the Harish-Chandra Schwartz space, and
regular semisimple orbital integrals of such functions converge
absolutely. The matrix entries of $A_{K_1,K_2}$ are such coefficients.
The open-correspondence identification underlying \refl{conv} is
$Z^\circ\simeq G\times G/T$, where
$p(g,aT)=g$ and $q(g,aT)=aya^{-1}g$; the relative density is the
one defining $\mu_y$. Consequently, for any nonnegative
$\varphi\in\calS(G)$ and any matrix entry, the absolute integral
against $p^*\varphi$ is
\eq{orbital-absolute}
 \int_{G(F)}\varphi(g)
       \int_{\Ome_y}
       |[A_{K_1,K_2}(xg)]_{ij}|\,d\mu_y(x)\,dg<\infty.
\end{equation}
Here and below a function on $G(F)$ is viewed as a half-density using
$|dg|^{1/2}$. To see the finiteness in \refe{orbital-absolute}, note
that $\pi(g)E_{K_2}$ takes only finitely many values as $g$ ranges
over $\operatorname{supp}(\varphi)$, since it is constant on right
$K_2$-cosets.

This proves the local integrability claim. More explicitly, for any
compact subset $C\subset Z(F)$, choose $\varphi\geq0$ which is positive
on $p(C)$. Under \refe{hecke-relative-line}, $p^*\varphi$ is a
nonvanishing positive test section near $C$, and every test section
supported in $C$ is bounded in absolute value by a constant multiple
of it. The complement $Z\setminus Z^\circ$ has measure zero, so the
operator-valued half-density defines a locally integrable
distribution on $Z$ with values in
$\Hom_{\CC}(V^{K_2},V^{K_1})$.

These distributions descend along $q$. Indeed, the two pullbacks to
the smooth variety $Z\times_{\calU}Z$ agree on the dense open subset
lying above $G$. Both are locally integrable, since smooth pullback
preserves local integrability. Their complement has measure zero,
so the two distributions agree everywhere. Smooth descent therefore
gives a linear map
\[
 \Psi_{K_1,K_2}:\calS_{1/2}(\calU)
 \longrightarrow\Hom_{\CC}(V^{K_2},V^{K_1}).
\]
Equivalently, one may use an admissible presentation of $\calU$:
the descended distributions are invariant and hence define
maps on the coinvariants defining its Schwartz space.

The maps $\Psi_{K_1,K_2}$ are compatible with averaging and with
the $G(F)\times G(F)$-action. These compatibilities hold on $G$ and
hence hold for their locally integrable extensions. If
$s\in\calS_{1/2}(\calU)$ is fixed by $J\times J$, its blocks are
therefore those of the single operator $\Psi_{J,J}(s)$, extended
from $V^J$ to $V$ by $E_J$. Explicitly, with this extension understood,
\[
 \Psi_{K_1,K_2}(s)=E_{K_1}\Psi_{J,J}(s)E_{K_2}.
\]
In particular, these blocks define an actual element of $\End'(V)$,
independent of $J$. Since every $s$ is fixed by some $J\times J$,
we obtain a $G(F)\times G(F)$-equivariant map
\[
 \Psi_\pi:\calS_{1/2}(\calU)\longrightarrow\End'(V)
\]
extending $f\mapsto\pi(f)$.

The quotient $\calS_{1/2}(\calU)/\calS(G)$ has a finite filtration
whose successive quotients are the nontrivial stratum modules
$\calS_{1/2}(X_\lambda)$ occurring in $\calU$. By \reft{exti},
\[
 \Hom_{G(F)\times G(F)}
 \bigl(\calS_{1/2}(\calU)/\calS(G),\End'(V)\bigr)=0.
\]
Thus the extension to $\calS_{1/2}(\calU)$ is unique, and
$\Psi_\pi=\Phi_\pi|_{\calS_{1/2}(\calU)}$.

Choose $K$ such that $f$ is bi-$K$-invariant. Then $h_D(f)$ is
$K\times K$-invariant, so its image is determined by its block on
$V^K$. Since $h_D(f)=q_!p^*f$, the construction above gives
\[
 \Phi_\pi(h_D(f))
 =\int_{G(F)}f(g)\int_{\Ome_y}
       E_K\pi(xg)E_K\,d\mu_y(x)\,dg
 =\pi(\mu_y)\pi(f).
\]
Here all integrals take values in $\End_{\CC}(V^K)$ and converge
absolutely by \refe{orbital-absolute};
the resulting operators are extended to $V$ using $E_K$. The last
equality follows from the definition of $\pi(\mu_y)$ and the
bi-$K$-invariance of $f$. This proves \refe{orbital-comparison}.
\epr

We now deduce \reft{main}(2) from \refp{orbital-comparison}. Choose
a compact open subgroup $K\subset G(F)$ with $V^K\neq0$.
Then by we have 
\refe{orbital-comparison},
\[
 \eta_D(\pi)\pi(e_K)
 =\Phi_\pi(h_D(e_K))=\pi(\mu_y)\pi(e_K).
\]
Now, for a square integrable $\pi$ we have
\[
 \pi\bigl((p_y)_*(dg)\bigr)
 =\frac{\chi_\pi(y)}{d_\pi(dg)}\operatorname{Id}_V
\]
(cf. Section 1 of \cite{Arthur} and Corollary 5 in \cite{RS}), where orbital integration is
defined by finite-dimensional averaging as above. Since
$\mu_y=c(y,dg)(p_y)_*(dg)$, this identity and
the preceding equality give
\[
 \eta_D(\pi)\pi(e_K)
 =c(y,dg)\frac{\chi_\pi(y)}{d_\pi(dg)}\pi(e_K).
\]
As $\pi(e_K)\neq0$, we obtain
$\eta_D(\pi)=c(y,dg)\chi_\pi(y)/d_\pi(dg)$.

\sec{}{Proof of \reft{cy}}

\ssec{}{The Hecke measure and the trace pairing}
We are going to compute the measure on $\Ome_y$ coming from the definition of Hecke operators explicitly.
Let
\[
 q_y:G/T\longrightarrow G,\qquad q_y(gT)=gyg^{-1}.
\]

Embed $E$ into $\End_F(E)$ by multiplication. The restriction of
the matrix trace pairing $B(X,Z)=\Tr(XZ)$ to $E$ is the field trace pairing.
It is nondegenerate since $E/F$ is separable. Thus
$$
 \End_F(E)=E\oplus E^\perp,\qquad
 V:=E^\perp\simeq \End_F(E)/E\simeq\Lie(G)/\Lie(T).
$$
Let $\nu_E$ be the invariant density on $G/T$ whose value at the identity
coset, in a tangent basis $v_1,\ldots,v_D$, where $D=N(N-1)$, is
\begin{equation}\label{eq:trace-density}
 \left|\det\bigl(\Tr(v_i v_j)\bigr)_{i,j=1}^{D}\right|_F^{1/2}.
\end{equation}

\th{}
\begin{equation}\label{eq:intrinsic}
 \mu_y=|\mathcal D_G(y)|_F^{1/2}(q_y)_*\nu_E.
\end{equation}
\eth
\prf
The part of the Hecke correspondence where both projective bundles are
trivial is parametrized by a point of $G$ and a conjugate of the
multiplication operator $y$, with the latter parameter retained in $G/T$.
Indeed, using vector bundles, the modification has the form
\[
 \calO^N\longrightarrow\calO(1)^N,\qquad A(t)=tI-C,
\]
after normalizing its invertible leading coefficient; its determinant
prescribes the characteristic polynomial of $C$ to be the minimal
polynomial of $y$. All such $C$ form the homogeneous space $G/T$.
Its two endpoint values give the two framed projective bundles.

Compute the canonical-line comparison after passing to a splitting field,
where
\[
 A(t)=\diag(t-y_1,\ldots,t-y_N),\qquad y_i\ne y_j\quad(i\ne j).
\]
This is an algebraic determinant calculation; its result will subsequently
be descended to $F$, before taking absolute values.
In the common meromorphic endomorphism bundle put
\[
 \mathcal K=\End(\mathcal E_1)\cap
 A^{-1}\End(\mathcal E_2)A,
 \qquad \mathcal E_1=\calO^N,\quad\mathcal E_2=\calO(1)^N.
\]
Twist all these sheaves by $\calO(-0-\infty)$. Since the two endpoint
endomorphism bundles are trivial and have no framed infinitesimal
automorphisms, the quotient exact sequences identify the two relative
tangent spaces $R_1,R_2$ with the respective spaces of local quotients by
$\mathcal K$. Scalar directions occur identically on both sides and cancel.

We use the canonical-line isomorphism of \cite[Theorem 1.1 and
Section 5.4]{EFK}. Its square compares the determinants of the two
deformation complexes. In the present model, the common terms cancel
in the two quotient exact sequences above, leaving the determinants
of the local quotients. The trace pairing and the residue description
of Serre duality identify these quotient determinants. Thus the squared
canonical-line isomorphism is computed by the determinant of their
local residue pairing. To check its normalization, in the local minuscule
model $A(z)=\diag(z,1,\ldots,1)$ the paired quotient generators are
$E_{i1}$ and $E_{1i}$, and
\[
 \Res_{z=0}\Tr\bigl(E_{i1} A(z)^{-1}E_{1i} A(z)\bigr)\,dz=1.
\]
The matrix units are dual for the trace pairing, and the displayed
residue is the Serre-duality evaluation in these local generators.
Taking determinants therefore introduces no additional scalar; taking
the absolute square root gives the half-density isomorphism. The
cotangent-line factor is the negative power specified in Section 1.2;
the chosen differential trivializes its density line and fixes the
numerical normalization used below.

For the matrix $A(t)$ above, the first quotient in the $ij$-direction is
supported at $y_j$. It pairs with the $ji$-direction in the second quotient.
Using $dt/t$, their pairing has coefficient
\begin{align}
 \beta(E_{ij},E_{ji})
 &=\Res_{t=y_j}
    \Tr\bigl(E_{ij}A(t)^{-1}E_{ji}A(t)\bigr)\frac{dt}{t}\notag\\
 &=\Res_{t=y_j}\frac{t-y_i}{t-y_j}\frac{dt}{t}
   =\frac{y_j-y_i}{y_j}.\label{eq:residue}
\end{align}
In this expression the residue evaluates the tangent-line-valued local
pairing against $dt/t$ at the modification point.

For completeness, the base tangent determinant introduces no extra
factor. Using the two framing actions, the tangent space to the
correspondence at this point is
$(\mathfrak g\oplus\mathfrak g)/\diag\mathfrak t$, and the two
projection differentials may be written
\[
 (u,v)\longmapsto v-u,
 \qquad (u,v)\longmapsto v-\Ad(y)u.
\]
In a root direction on which $\Ad(y)$ acts by $a$, bases for their kernels
are $(1,1)$ and $(a^{-1},1)$, respectively; the vector $(0,1)$ is a lift
of a base tangent generator for both projections. The comparison of
these two exact-sequence determinants therefore contributes a power
$a^{\pm1}$. The product over all root directions is $1$, because
$\det(\Ad(y)\mid\mathfrak g/\mathfrak t)=1$. On the torus directions
the two projections coincide.

Relative to the trace pairing on $V$, the determinant of
\eqref{eq:residue} is consequently
\begin{equation}\label{eq:weyl}
 \prod_{i\ne j}\frac{y_j-y_i}{y_j}
 =\prod_{i\ne j}(1-y_i/y_j)=\mathcal D_G(y).
\end{equation}
One can also see the direction of this determinant factor directly from
the tangent exact sequences. They identify the ratio of the two base
canonical lines with $\det R_1/\det R_2$. The residue pairing identifies
$\det R_1$ with $(\det R_2)^{-1}$, with coefficient
\eqref{eq:weyl}; hence the coefficient in the relative density is its
absolute square root, rather than its inverse.

The quotient sheaves, trace-residue pairing, and determinant
identifications are defined over $F$. Descending this calculation and
taking absolute square roots therefore gives
$|\mathcal D_G(y)|_F^{1/2}\nu_E$ on the parameter space $G/T$.
This proves \eqref{eq:intrinsic}.
\epr
\ssec{}{Comparison with Haar measures}
\prf
On $\GL(N,F)$ the invariant trace-pairing density is
\[
 dm_{\GL(N)}(g)=\frac{\prod_{i,j}dg_{ij}}{|\det g|_F^N},
\]
where additive measure on each matrix entry gives $\calO_F$ volume $1$.
Indeed, the trace pairing in the matrix-unit basis has determinant $\pm1$.
On $E^*$ the corresponding density is
\[
 dm_{E^*}(a)=\delta_E^{1/2}
                  \frac{da}{|N_{E/F}(a)|_F},
 \qquad da(\calO_E)=1.
\]
The factor $\delta_E^{1/2}$ is the square root of the absolute determinant
of the field trace pairing in an integral basis. Orthogonality of
$M_N(F)=E\oplus E^\perp$ shows that
\[
 \nu_E=\frac{dm_{\GL(N)}}{dm_{E^*}}
 \quad\text{on}\quad\GL(N,F)/E^*\simeq G(F)/T(F).
\]

Divide both group measures by $dm_{F^*}=dx/|x|_F$. This produces
measures $dg_{\mathrm{can}}$ on $G(F)$ and $dt_{\mathrm{can}}$ on $T(F)$, with
\[
 \nu_E=dg_{\mathrm{can}}/dt_{\mathrm{can}}.
\]
Counting invertible matrices modulo the maximal ideal gives
\[
 dm_{\GL(N)}(\GL(N,\calO_F))=\prod_{j=1}^{N}(1-q^{-j}),
 \qquad dm_{F^*}(\calO_F^*)=1-q^{-1}.
\]
Thus
\begin{equation}\label{eq:group-volume}
 dg_{\mathrm{can}}(K)=\prod_{j=2}^{N}(1-q^{-j}).
\end{equation}
Likewise
\[
 dm_{E^*}(\calO_E^*)=\delta_E^{1/2}(1-q^{-f}).
\]
The quotient $E^*/(F^*\calO_E^*)$ has order $e$.
Its cosets have equal quotient measure, so
\begin{equation}\label{eq:torus-volume}
 dt_{\mathrm{can}}(T(F))=
 e\,\delta_E^{1/2}\frac{1-q^{-f}}{1-q^{-1}}.
\end{equation}
Let $dt_T$ be probability Haar measure on the compact group $T(F)$.
Equations \eqref{eq:group-volume} and \eqref{eq:torus-volume} imply,
for arbitrary $dg$,
\[
 \nu_E=
 \frac{dg_{\mathrm{can}}(K)}{dt_{\mathrm{can}}(T(F))\,dg(K)}\frac{dg}{dt_T}
 =\frac{A_E}{dg(K)}\frac{dg}{dt_T}.
\]
Finally, quotient integration and $dt_T(T(F))=1$ give
\[
 (q_y)_*(dg/dt_T)=(p_y)_*dg.
\]
Combining this identity with \eqref{eq:intrinsic} proves
\reft{cy}.
\end{proof}

\end{document}